\documentclass[letter,11pt]{article}
\usepackage{amsfonts}
\usepackage{amssymb}
\usepackage{amsmath}
\usepackage{cite}
\newtheorem{theo}{Theorem}
\newtheorem{prop}[theo]{Proposition}
\newtheorem{lema}[theo]{Lemma}

\newtheorem{obs}[theo]{Remark}
\newcommand{\R}{\mathbb{R}}
\newcommand{\D}{\mathbb{D}}
\newcommand{\T}{\mathbb{T}^1}

\newcommand{\C}{\mathbb{C}}
\newcommand{\N}{\mathbb{N}}
\newcommand{\eps}{\varepsilon}
\newcommand{\Z}{\mathbb{Z}}

\title{Towards a semilocal study of parabolic invariant curves for fibred holomorphic maps}
\author{Mario PONCE\\ PUC-Chile}
\begin{document}
\maketitle
\begin{abstract}
We introduce the study of the local dynamics around a parabolic indifferent invariant curve for fibred holomorphic maps. As in the classical non-fibred case, we show that petals are the main ingredient. Nevertheless, one expects the properties of the base rotation number should play an important  role in the arrangement of the petals. We exhibit examples where the existence and the number of petals depend not just on the complex coordinate of the map, but on the base rotation number. Furthermore, under additional hypothesis on the arithmetics and smoothness of the map, we present a  theorem that allows to characterize the local dynamics around a parabolic invariant curve.
\end{abstract}
\section{Introduction}
Among the first results on one-dimensional complex dynamics one find the understood of the  behavior of points near a fixed point under iteration of a holomorphic map.  Since the works of K\"onig and Poincar\'e, one knows that the derivative of the map at the fixed point rules the local dynamics in the hyperbolic case. When the derivative at the fixed point is not hyperbolic, one says that the fixed point is indifferent and finer techniques are needed, because the very fast exponential convergence due to hyperbolicity is no longer available. 
\\

When the derivative is a    rotation by a rational angle,  one says that the fixed point is rationally indifferent. An elegant result due to Leau and Fatou states that the local dynamics is actually ruled by the derivative at the fixed point. Indeed, the local dynamics follows a nice and simple pattern known as a Leau-Fatou flower:

 Let $f(z)=z+az^{n+1}+\dots, \ a\neq 0$. A simply connected open set $U$ is called an attracting petal for $f$ if \[
 0\in \partial U\ , \ f(\overline U)\subset U\cup\{0\}\ \textrm{and} \ \bigcap_{j\geq 0}f^j(\overline U)=\{0\}.
 \] 
One defines the notion of a repelling petal in an analogous way. 
 \begin{theo}[Leau-Fatou flower, see \cite{LEAU97, FATO19, MILN06, CAGA93}] Let $f $ be as above. There exists $n$ disjoint attracting petals and $n$ disjoint repelling petals. These $2n$ petals, together with the fixed origin itself, form a neighborhood of the origin. 
 \end{theo}
Recently, F. Le Roux \cite {LERO08} has shown a topological version of this result for surface homeomorphisms with Lefschetz index bigger than $1$. \\

In the irrational indifferent case, the local dynamics is ruled by the derivative for the most cases. In fact, by assuming an extra arithmetical hypothesis on the rotation angle,  the local dynamics is conjugated to the corresponding rigid rotation  (the Brjuno arithmetical condition has full Lebesgue measure, see for instance \cite{YOCC95}).  
\\

In the last years, a special interest have been paid to skew-product dynamical systems, and many results about complex skew-product dynamics have appeared (see \cite{JONS99}, \cite {SEST99}, \cite{SUMI00}). On the one hand, skew-product dynamics presents a rich source of examples, and in the other hand, it is considered as an intermediate step towards the higher dimensional complex dynamics. This paper is devoted to the study of continuous maps in the form
\begin{eqnarray*}
F:\T\times \D&\longrightarrow& \T\times \C\\
(\theta, z)&\longmapsto&\left(\theta+\alpha, f_{\theta}(z)\right)
\end{eqnarray*}   
where $\T=\R/\Z$, $\alpha$ is an irrational angle and each $f_{\theta}:\D\to \C$ is holomorphic. We call these as {\it fibred holomorphic maps} (abbreviated as {\it fhm}). Since $\theta\mapsto \theta+\alpha$ is minimal, $F$ has no fixed points. However, the  natural object  playing the role of a fixed point is an {\it invariant curve}, that is, a continuous curve $u:\T\to \D$ such that 
\[
F\left(\theta, u(\theta)\right)=\left(\theta+\alpha, u(\theta+\alpha)\right).
\]
The aim of this  paper is to understand orbits of points close to an invariant curve. In \cite{PONC07}, the author shows that a Poincar\'e-K\"onig Theorem holds for transversally hyperbolic invariant curves. Also, a version of the Siegel's Theorem holds for Diophantine transversal rotation number. \\

In this work we introduce the study of the local dynamics around a parabolic indifferent invariant curve. As in the classical non-fibred case, we find out that petals are the main ingredient (see Proposition \ref{flor}). Nevertheless, one expects the properties of the base rotation number should play an important  role in the arrangement of the petals. We exhibit examples where the existence and the number of petals depend not just on the complex coordinate of the map, but on the base rotation number (see Section \ref{ejemplos}). Furthermore, under additional hypothesis on the arithmetics and smoothness of the map, we present a  theorem that allows us to characterize the local dynamics around a parabolic invariant curve (see Theorem \ref{identidaddiofantina}).\\

Given a curve $u:\T\to \D$ which is invariant for $F$, we can perform a fibred translation in order to locate $u$ at the zero section $\T\times \{0\}$. Therefore, in the sequel we will deal only with complex maps in the form $f_{\theta}(z)=a_1(\theta)z+a_2(\theta)z^2+\dots$. The invariant curve will be called {\it parabolic indifferent} if $a_1(\theta)\equiv \lambda$, with $\lambda $ being a constant primitive root of the unity. 
\begin{obs}
Under additional hypothesis on the arithmetics and smoothness of the map $F$, we can relax the above definition just by demanding the number
\[
\frac{1}{2\pi i }\int \log a_{1}(\theta) d\theta
\] to be a rational number. Indeed, by solving a cohomological equation and performing a fibred homotopy we can recover the above form for the map $f_{\theta}$. A natural direction of research is to determine whether or no the results contained in this work still hold for this weak definition of a parabolic curve when those additional hypothesis are absent.
\end{obs} 
\noindent{\it Acknowledgments. } I wish express my gratitude to  Jan Kiwi and Jean-Christophe Yoccoz for many helpful discussions and comments. Part of this work was prepared thanks to the hospitality of   the Coll\`ege de France. This work was partially supported by FONDECYT 11090003 grant and MathAMSUD Project DySET.

\section{Fatou coordinates and cohomological equation}
In the one-dimensional case the study of the local complex dynamics near a parabolic fixed point is achieved by looking the map as acting on a neighborhood of infinity. More precisely, one introduces a conjugacy by the map $i(z)=-z^{-1}$, the so called Fatou coordinates. This procedure will also give us  good information about the local dynamics of a fibred holomorphic map on a neighborhood of a parabolic invariant curve. Let's consider the fibred Fatou coordinates
\begin{eqnarray*}
I:\T\times \overline{\C}&\longrightarrow&\T\times \overline{\C}\\
(\theta, z)&\longmapsto&(\theta, Z=-z^{-1}).
\end{eqnarray*}
Let $F$ be a {\it fhm} in the form
 \begin{equation}\label{dinamica}
 F(\theta,z)=\left( \theta+\alpha, z+a_2(\theta)z^2+a_3(\theta)z^3+\dots\right )
 \end{equation}
defined on $\T\times D(0,r)$ for some $r>0$. By looking at infinity one gets
 the fibred map
 $
  \tilde F=I\circ F \circ I^{-1}
  $
  which is defined on  $\T\times I(D(0,r))$ and takes the form
  \begin{equation*}
  \tilde F(\theta,Z)=\left( \theta+\alpha, \frac{-1}{\frac{-1}{Z}+\frac{a_2(\theta)}{Z^2}-\frac{a_3(\theta)}{Z^3}+\dots}\right)
  =\left( \theta+\alpha, Z\left[\frac{1}{1-\frac{a_2(\theta)}{Z}+\frac{a_3(\theta)}{Z^2}-\dots}\right]\right).
  \end{equation*}
As the convergence radius of the series $a_2(\theta)z^2+a_3(\theta)z^3+\dots$ is at least $r$ for every $\theta\in \T$, there exists $C>0$ such that 
\[
\left |\frac{a_2(\theta)}{Z}-\frac{a_3(\theta)}{Z^2}+\dots\right|<1 
\]
for $|Z|>C$. Thus, for $|Z|>C$ we have 
 \[
  \tilde F(\theta,Z)=\left( \theta+\alpha, Z+a_2(\theta)+\frac{a_2(\theta)^2-a_3(\theta)}{Z}+\dots\right).
  \]
  At a first glance, for  big enough $|Z|$ this dynamics is pretty like a fibred translation by $a_2(\theta)$. In fact, we will see that under some mild hypothesis, the map $\tilde F$ can actually be conjugated to a fibred translation. Let's explain the general procedure: we look for a conjugacy 
  \[
  T_c(\theta,Z)=(\theta,Z+c(\theta))
  \]
   in such a way that  $T_c\circ \tilde F\circ T_c^{-1}$ takes the form  
  \begin{equation}\label{primeraconjugacion}
  (\theta,Z)\longmapsto\left(\theta+\alpha,Z+k+\frac{b_{1}(\theta)}{Z}+\frac{b_{2}(\theta)}{Z^2}+\dots\right), 
  \end{equation}
  where $k=\int_{\T}a_2(\theta)d\theta$. This drive us to consider the cohomological equation
  \begin{equation}\label{cohomo}
  c(\theta+\alpha)-c(\theta)=-a_2(\theta)+k.
  \end{equation}
  This equation has been widely studied. It is known that a continuous solution for (\ref{cohomo}) exists if and only if the Birkhoff sums of the function $a_2(\cdot)-k$ are bounded (Gottschalk-Hedlund's Lemma, see \cite{GOHE55}). Moreover, such a solution is unique if we require $\int_{\T}c(\theta)d\theta=0$ (which will be always our choice). We don't pretend to discuss the precise and optimal hypothesis for the existence of continuous solutions for this equation. We just mention the  following result:
  \begin{lema}[see \cite{KAHA95}]\label{babykam}
  Let $\alpha$ be a $CD(\tau)$ Diophantine number, $a_2(\theta)$ a differentiable $C^{1+\tau}$ function. Equation (\ref{cohomo}) has a continuos solution $c:\T\to \C\quad_{\blacksquare}$
  \end{lema}
  For more precise results, the reader can refer to   Herman's works \cite{HERM79, HERM83-as1}. The next lemma is trivial, but useful in the sequel:
  \begin{lema}
  Let $\alpha$ be irrational and  let $a_2(\theta)$ be  a trigonometric polynomial.  Equation (\ref{cohomo}) has a  solution  which is a trigonometric polynomial $\quad_{\blacksquare}$
  \end{lema}
  \section{Fibred flowers}
  In this section we show that the classical theory of the Leau-Fatou flower also hold in the fibred setting. However, we need to know that the leading  coefficient has non-zero mean. In the sequel, we write $\int_{\T}\cdot$ instead of $\int_{\T}\cdot \ d\theta$.   
    \subsection{Case $z+Az^2+\dots$}
 We continue with the discussion on the previous section. Assume that  $k=\int_{\T}a_2\neq 0$. Equation (\ref{cohomo}) could provide the required fibred translation $T_c$. Of course, there are two situations we need to treat:  
  \paragraph{Equation (\ref{cohomo}) has continuous solution. } In this case, and 
  by conjugating by an extra homothetic map $A_k(\theta, Z)=(\theta, k^{-1}Z)$, the map $\tilde F$ can be thought as having the form 
  \begin{equation}\label{tildefuno}
  \tilde{F}(\theta,Z)=\left(\theta+\alpha, Z+1+\frac{b_{1}(\theta)}{Z}+\frac{b_{2}(\theta)}{Z^2}+\dots\right)
  \end{equation}
  defined and convergent for $|Z|>C>0$, where $C$ is a positive constant (possibly greater than $C$ defined above) and $b_{j}:\T\to \C$ are continuous functions. 
  \paragraph{Equation (\ref{cohomo}) has no continuous solution. }    Let $l:\T\to \C$ be  a trigonometric polynomial, such that
  \begin{itemize}
  \item[$i) $] $\hat l(0)=\int_{\T}l(\theta)=k$.
  \item[$ii) $] $|l(\theta)-a_2(\theta)|<\frac{k}{1000}$ for every $\theta \in \T$.
  \end{itemize}
    We solve the cohomological equation
    \[
    c(\theta+\alpha)-c(\theta)=-l(\theta)+k
    \]
  and obtain a trigonometric polynomial $c:\T\to \C$ which provides a conjugacy $T_c$ such that
  \[
  T_c\circ \tilde F\circ T_c^{-1}(\theta, Z)=\left(\theta+\alpha, Z+k+\left(l(\theta)-a_2(\theta)\right)+\frac{b_1(\theta)}{Z}+\dots\right).
  \]
  An extra homothetic conjugacy allows us to consider $\tilde F$ having the form
  \[
  \tilde F(\theta, Z)=\left(\theta+\alpha, Z+1+\tilde l(\theta)+\frac{b_1(\theta)}{Z}+\frac{b_2(\theta)}{Z^2}+\dots\right)
  \]
  where $|\tilde l(\theta)|<<1$ and $b_{j}:\T\to \C$ are continuous functions.\\

  In any case, we get a form for the fibred map that allows us to follow the classical construction providing the attracting and repelling petals. We perform the fibred version of this  procedure in the next section.
  \paragraph{Near infinity dynamics and invariant regions. }
   Let $A>0$ be a constant such that the regions
   \[
   \Omega_A^+=\big\{Z=x+iy\big| x>A-|y|\big\}\quad ,\quad \Omega_A^-=\big\{Z=x+iy\big| x<-A+|y|\big\}
   \]
   are contained in $\C\setminus D(0,C)$. The region $\Omega_A=\Omega_A^+\cup \Omega_A^-$ is a neighborhood of $\infty$. For $C_2$ large enough and $Re(Z)>C_2$ we have
   \[
   |Re\tilde F^n_{\theta}(Z)|>Re(Z)+\frac{n}{2}.
   \]
   Note that $\tilde F(\Omega_A^+)\subset \Omega_A^+$. 
   \\
   
   Now, let's consider the region $\mathcal R_{\theta}\subset \C$, delimited by the line $\mathcal{L}=\{x=L\}$, with $L>A$, and its image $\mathcal{L}_{\theta}=\tilde F_{\theta-\alpha}(\mathcal{L})$. For each $\theta\in \T$ there exists a homeomorphism $H_{\theta}$ from $\mathcal{R}_{\theta}$ into the region $\{Z\in  \C\ \big |\  0\leq Re(Z)\leq 1\}$. We can pick $H_{\theta}$ verifying $H_{\theta}(L+iy)=iy$ and $H_{\theta}\left(\tilde F_{\theta-\alpha}(L+iy)\right)=1+iy$.  Moreover, we can take the map $H(\theta,Z)=(\theta,H_{\theta}(Z))$ to be continuous. If we restrict the dynamics of $\tilde F$ to $\T\times \{Re(Z)\geq L\}$, the subset $\bigcup_{\theta\in \T}\{\theta\}\times \mathcal{R}_{\theta}$ becomes a fundamental domain and we can extend the homeomorphisms $H_{\theta}$ in order to get a homeomorphism
   \begin{eqnarray*}
   H:\T\times \{Re(Z)\geq L\}&\longrightarrow&\T\times\{Re(Z)\geq 0\}\\
   (\theta,Z)&\longmapsto&\left(\theta,H_{\theta}(Z)\right)
   \end{eqnarray*}
   verifying 
   \[
   H\circ \tilde F\circ H^{-1}(\theta,Z)=(\theta+\alpha, Z+1)
   \]
   for every $(\theta,Z)\in \T\times \{Re(Z)\geq 0\}$. Indeed, for $(\theta,Z)$ with $Re(Z)\geq 0$, there exists a unique $n\in \N$ such that $0\leq Re(Z-n)<1$. Then we define 
   \[
   H^{-1}(\theta,Z)=\left(\theta, \tilde{F}^n_{\theta-n\alpha}\left(H_{\theta-n\alpha}^{-1}(Z-n)\right)\right).
   \]
   By using an analogous construction we can conjugate $\tilde{F}^{-1}$ to $(\theta-\alpha, Z-1)$ into $Re(Z)<-L$.
   
    If we look at the regions $\Omega_A^+, \Omega_A^-$ in the $z-$plane (the original coordinates in a neighborhood of the origin) we get, for each fibre,  two topological discs $\mathcal{P}_{\theta}^+, \mathcal{P}_{\theta}^-$ depending continuously on $\theta$, such that
    \begin{enumerate}
    \item $0\in \partial \mathcal{P}_{\theta}^+$.
    \item $\mathcal{P}_{\theta}^+\cup \mathcal{P}_{\theta}^-$ is a neighborhood of the origin.
    \item $\tilde F_{\theta}(\mathcal{P}_{\theta}^+)\subset \mathcal{P}_{\theta+\alpha}^+$ and $\tilde F^{-1}_{\theta}(\mathcal{P}_{\theta}^-)\subset \mathcal{P}_{\theta-\alpha}^+$.
    \end{enumerate}
    The disjoint union of the sets $\{\theta\}\times \mathcal{P}_{\theta}^+$ forms an open tube having the invariant zero section in the boundary. Moreover, this tube is  forward invariant and converges to the invariant curve by forward iteration. We call such a set a {\it fibred attracting petal} (we define analogously the notion of a  {\it fibred repelling petal}).
    Summarizing, we can state the 
    \begin{prop}\label{propunpetalo}
 Let $F$ be a fibred holomorphic map as in (\ref{dinamica}). If $\int_{\T}a_2\neq 0$  then there exists one fibred attracting petal and one fibred repelling petal. The union of this two sets form a tubular neighborhood of the invariant curve $\quad_{\blacksquare}$
    \end{prop}
    \paragraph{Exterior direction of fibred petals. }
   We have shown that the original map $F$ is conjugated to the map $(\theta, z)\mapsto\left(\theta+\alpha, z+z^2+\dots \right) $ via the fibred homeomorphism $L=A_k\circ T_c\circ I$. Explicitly one has 
   \[
   L(\theta, z)=\left(\theta, k^{-1}\left(\frac{-1}{z}+c(\theta)\right)\right).
   \]
   In Section \ref{an} we will need to cut the plane $\C$ along the repulsive direction of the fibred petal, hence we compute it here: In the coordinates of $Z$ this direction corresponds to the direction of the curve $\gamma_t:t\in \R\mapsto t\in \C$ with $t \to +\infty$. Thus, the direction we look for is the direction of $L^{-1}(\theta, \gamma_t)$ as $t\to +\infty$: 
   \begin{eqnarray*}
   \frac{d}{d t} \frac{-1}{kt-c(\theta)}&=&\frac{k}{(kt-c(\theta))^2}\\
   &=&\frac{1}{kt^2}\frac{1}{\left(1+\frac{c(\theta)}{kt}\right)^2}
   \end{eqnarray*}
   that is parallel to $k^{-1}$. Thus, the exterior direction of the fibred petal does not depend on $\theta$ and corresponds to the direction of $k^{-1}$. 
   \subsection{Case $z+Az^{n+1}+\dots$}\label{an}
   In this section we consider the following fibred map
 \begin{equation}\label{dinamica_n}
 F(\theta,z)=\left( \theta+\alpha, z+a_{n+1}(\theta)z^{n+1}+a_{n+2}(\theta)z^{n+2}+\dots\right )
 \end{equation}
 with $n> 1$.  We assume  that $\int_{\T}a_{n+1}\neq 0$. Write  
 $
 \int_{\T}a_{n+1}=re^{i\Theta}
 $
 with $r>0$ and $\Theta\in [0, 2\pi)$. We consider the {\it repulsive directions} as the $n$  unitary vectors $e_j=e^{i\left(\frac{2\pi j}{n}-\Theta\right)}$ for $0\leq j\leq n-1$.  We call $\mathcal{R}_j$ the open region limited by two consecutive repulsive directions $e_j, e_{j+1}$. Pick $\mathcal{R}=\mathcal{R}_j$ to be one of these regions. We define the homeomorphism
 \begin{eqnarray*}
 \Phi_n:\T\times \mathcal{R}&\longrightarrow& \T\times \{\C\setminus \mathcal{L}_{-\Theta}\}\\
 (\theta, z)&\longmapsto& (\theta,w= z^n)
 \end{eqnarray*}
  where $\mathcal{L}_{-\Theta}=\R_+ e^{-i\Theta}$. We conjugate $F$ by $\Phi_n$   obtaining a fibred holomorphic map on $\T\times \{\C\setminus \mathcal{L}_{-\Theta}\}$ 
  \[
  (\theta, w)\longmapsto\left(\theta+\alpha, w+na_{n+1}(\theta)w^2+\dots\right).
  \]
  Since $\int_{\T}na_{n+1}\neq 0$ Proposition \ref{propunpetalo} implies this map presents one fibred attracting petal and one fibred repelling petal (cutted by the line $\mathcal{L}_{-\Theta}$). Pasting together the regions $\mathcal{R}_j$ we get the
  \begin{prop}\label{flor}
  Let  $F$ be a fibred holomorphic map as in (\ref{dinamica_n}). If $\int_{\T}a_{n+1}\neq 0$ then there exists $n$ fibred attracting petals and $n$ fibred reppelling petals. The union of these $2n$ fibred petals and the invariant curve form a tubular neighborhood of the invariant curve$\quad_{\blacksquare}$
  \end{prop}
  \subsection{Case $\lambda z+Az^2+\dots, 	\ \lambda^n=1$}
  Let 
   \begin{equation*}
 F(\theta,z)=\left( \theta+\alpha, \lambda z+a_{2}(\theta)z^{2}+a_{3}(\theta)z^{3}+\dots\right )
 \end{equation*}
 where $\lambda$ is a $n^{th}$ primitive root of the unity. The iterate $F^n$ of $F$ has the form
  \begin{equation*}
 F^n(\theta,z)=\left( \theta+n\alpha, z+b_{2}(\theta)z^{2}+b_{3}(\theta)z^{3}+\dots\right ),
 \end{equation*}
hence previous results (and the procedure at Section \ref{infinitoreducible}) can be applied. In the very same way as in the classical non fibred case, if $F^n$ presents a $p$ petals parabolic behavior, let's say $\mathcal{P}_1, \mathcal{P}_2, \dots , \mathcal{P}_p$ are the fibred  attracting petals, then $F(\mathcal{P}_j)$ equals $\mathcal{P}_i$ for some $i$ (in the sense of germs). Thus, $F$ permutes the petals in cycles of length $n$. We conclude that $n$ divides $p$. 

%
    \section{Case $\int_{\T}a_2=0$}\label{inta20}
    
   Let's come back to the situation of the map (\ref{primeraconjugacion}) and assume   that $\int_{\T}a_2=0$.  Suppose that a continuous solution $c_2:\T\to \C$ exists for the cohomological equation (\ref{cohomo}). By performing the fibred translation and then coming back to the original coordinates, it is not hard to see that the map $F$ is conjugated to 
   \begin{equation*}
F_{(1)}(\theta,z)\longmapsto\left(\theta+\alpha,z-d_1z^3+d_2z^4+(d_1^2-d_3)z^5\dots\right)
   \end{equation*}
   where
   \begin{eqnarray*}
   b_1&=&a_2^2-a_3\\
   b_2&=&a_4-2a_2a_3+a_2^2\\
   b_3&=&-a_5+a_3^2+2a_2a_4-3a_2^2a_3+a_2^4\\
   d_1&=&b_1\\
   d_2&=&b_1c_2+b_2\\
   d_3&=&b_1c_2^2+2b_2c_2+b_3.
   \end{eqnarray*}
   We call  $F_{(1)}$ the {\it order $2$-reduction} of $F=F_{(0)}$.  In general, if the leading coefficient (let's say $a_{n+1}$) verifies $\int_{\T}a_{n+1}=0$, and the corresponding cohomological equation has a continuous solution,  we can perform a reduction in order to get a fibred map with leading coefficient of order $\tilde n> n+1$. We call it the order $n+1$-reduction (see Section \ref{infinitoreducible} for details). We can easily obtain formulas for the new coefficients $\{\tilde a_j\}_{j\geq n+2}$ in terms  of the old ones $\{a_j\}_{j\geq n+1}$ and the solution $h_{n+1}$ to the corresponding cohomological equation. For example:
      \begin{enumerate}
      \item[$n+2)$] Provided that $n\geq 2$ one has
   \[
   \tilde a_{n+2}=a_{n+2}.
   \]
   \item[$n+3)$] Provided that $n\geq 3$ one has
   \[
   \tilde a_{n+3}=a_{n+3}.
   \]
   In  the case $n=2$ one has:
   \begin{enumerate}
   \item[$2+3)$] \[
   \tilde a_{5}=3a_{3}h_{3}+a_{5}.
   \]
   \end{enumerate}
   \end{enumerate}

   \subsection{Some pathological examples}\label{ejemplos}
In this section we exhibit some examples showing new phenomena on parabolic fibred maps. First three examples show that counting petals is not easy {\it a priori}, in contrast with the well known parabolic  theory of the one-dimensional complex dynamics. Last two examples are toy models for more complicated behavior.
\begin{enumerate}
\item Let's take
\[
F(\theta, z)=\left(\theta+\alpha, z+\sin(\theta)z^2\right).
\]
In this case  $\int_{\T}a_2=0$ and the cohomological equation has continuous solution for every irrational $\alpha$. By performing the order $2$-reduction, we find out  that  $F$ is topologically equivalent to
\[
(\theta, z)\longmapsto \left(\theta+\alpha, z-\sin^2(\theta)z^3+\dots\right).
\]
Since $\int_{\T}\sin^2(\theta)d\theta\neq 0$, $F$ has a fibred flower with $2$ attracting petals and $2$ repelling petals . Also note  that every coefficient $a_j$ of $F$ verifies $\int_{\T}a_j=0$.
\item Let's take
\[
F(\theta, z)=\left(\theta+\alpha, z+\sin(\theta)z^2+\sin(\theta)^2z^3+\cos(\theta)^2z^4\right).
\]
Note $\int_{\T}a_2=0$ and $\int_{\T}a_3\neq 0$. We could  guess a $2$ petals parabolic behavior. Nevertheless, the order $2$-reduction says that $F$ is topologically equivalent to 
\[
(\theta, z)\longmapsto \left(\theta+\alpha, z+\left(1-2\sin^3(\theta)\right)z^4+\dots\right).
\]
Hence a $3$ petals parabolic behavior appears! 
\item In the next example, we will construct a parabolic fibred map presenting petals, but in such a way that the number of petals depends on the base rotation number. This suggests  that the number of petals can not be computed by means of an integral just depending on the complex coordinate of the map (as in the one-dimensional case). We use the notations of the beginning of Section \ref{inta20}.  We will construct a fibred polynomial of degree $5$ with trigonometric polynomial coefficients:
\[
F_{\alpha}(\theta,z)=\left( \theta+\alpha, z+a_2(\theta)z^2+a_3(\theta)z^3+a_4(\theta)z^4+a_5(\theta)z^5\right ).
\]
We use the sub-script  for making explicit the dependence on the base rotation number.
Take  a non identically zero $a_2$ such that $\int_{\T}a_2=0$.  Let $h_2^{\alpha}$ be the solution to the corresponding  cohomological equation. Note that $h_2^{\alpha}$ depends on $\alpha$ (we put $\alpha$ as a super-script in order to remark this).  Take $a_3$ being constant and equal to $\int_{\T}a_2^2$. The order $2$-reduction gives
\[
\tilde F_{\alpha}(\theta, z)=\left(\theta+\alpha, z-d_1z^3+d_2^{\alpha}z^4+(d_1^2-d_3^{\alpha})z^5+\dots\right).
\]
The above choice of $a_3$ yields $\int_{\T}d_1=0$. We can make the order $3$-reduction and construct the solution $h_3^{\alpha}$ to the corresponding cohomological equation.  We conjugate the original map to
\begin{equation}\label{oderhand}
(\theta, z)\longmapsto \left(\theta+\alpha, z+d_2^{\alpha}z^4+\left [-3d_1h_3^{\alpha}+d_1^2-d_3\right ]z^5+\dots\right).
\end{equation}
Recall \[d_2^{\alpha}=b_1h_2^{\alpha}+b_2=(a_2^2-a_3)h_2^{\alpha}+a_4-2a_2a_3+a_2^2.\]
Fix $\alpha=\alpha^*$ and pick $a_4$ to be a constant and such that $\int_{\T}d_2^{\alpha^*}=0$. By performing the order $4$-reduction we can conjugate $F_{\alpha^*}$ to
\[
(\theta, z)\longmapsto \left(\theta+\alpha^*, z+\left [-3d_1h_3^{\alpha^*}+d_1^2-d_3\right ]z^5+\dots\right).
\]
Finally, choose $a_5$ (and so $d_3$) such that $\int_{\T} [-3d_1h_3^{\alpha^*}+d_1^2-d_3 ]\neq 0$. Hence $F_{\alpha^*}$ presents a $4$ petals parabolic behavior.

By the other hand, pick $\alpha=\alpha^{**}$ such that in (\ref{oderhand}) we have $\int_{\T}d_2^{\alpha^{**}}\neq 0$. Then, the reduction procedure stops there and $F_{\alpha^{**}}$ presents a $3$ petals parabolic behavior.
\item  Let $\alpha$ be an irrational number not belonging to the Brjuno class (see \cite{YOCC95} for definitions). The following theorem is due to Yoccoz
\begin{theo}[see \cite{YOCC95}]
The quadratic polynomial $P(z)=e^{2\pi i  \alpha}z+z^2$ is not linearizable. Furthermore, there exists periodic orbits approximating the fixed point $z=0\quad_{\blacksquare}$
\end{theo}
Consider the following fibred holomorphic map:
\[
Q(\theta, z)=\left(\theta+\alpha, P(z)\right).
\]
Of course, this fibred map does not present a parabolic behavior since there are periodic curves converging to the invariant curve. By performing the (not isotopic to the identity) change of coordinates $(\theta, z)\mapsto (\theta, e^{-2\pi i  \theta}z)$ we get the map
	\[
	(\theta, z)\mapsto \left(\theta+\alpha, z+e^{2\pi i (\theta-\alpha)}z^2\right)
	\]
	which seems to be  a  fibred parabolic dynamics. Note that $\int_{\T}a_j=0$ for every $j\geq 2$. We will come back to this example in Section \ref{infinitoreducible}.
	\item The next one is an interesting example, since, at some extent, it should model the non reducible case, being the source of new and rich phenomena in fibred holomorphic maps. Let $a:\T\to \C$ be a continuous function and such that $\int_{\T}a=0$. We consider the fibred map defined in $\T\times \overline{\C}$
	\begin{eqnarray*}
	F_a(\theta, z)&=&\left(\theta+\alpha, \frac{z}{1-a(\theta)z}\right)\\
	&=&\left(\theta+\alpha, z+a(\theta)z^2+a(\theta)^2z^3+\dots\right).
	\end{eqnarray*}
	By the change of coordinates at infinity we get the map
	\begin{equation}\label{nosee}
	\tilde F_a(\theta, Z)=\left(\theta+\alpha, Z+a(\theta)\right).
	\end{equation}
	In the literature, (\ref{nosee}) is known as {\it cylindrical cascade} and has been widely studied (see for instance \cite{GOHE55}, \cite{BESI51},  \cite{ATKI78}, and the fairly complete  introduction on\cite{FRLE10}).
	\\
	
	If the corresponding cohomological equation 
	\begin{equation}\label{otracoho}
	c(\theta)-c(\theta+\alpha)=a(\theta)
	\end{equation}
	has a continuous solution, then the map $\tilde F_a$, and a posteriori $F_a$, is topologically equivalent to the fibred identity map $(\theta, Z )\mapsto (\theta+\alpha, Z)$.  On the other hand, let's concentrate in the situation where we can not solve the cohomological equation. In this case one says that the dynamical system (\ref{nosee}) is {\it non-integrable}. The most outstanding result concerning topological dynamics of cylindrical cascades is
	\begin{theo}[Atkinson '78, \cite{ATKI78}] Let $a:\T\to \C$ be a  continuous function with $\int_{\T}a=0$. If   $\tilde F_a$ is non-integrable then there exists a non-zero complex number $\tau $ such that the real cylindrical cascade \[(\theta, t)\longmapsto (\theta+\alpha, t+<\tau, a(\theta)>)\] is topologically transitive $\quad_{\blacksquare}$
	\end{theo}
	Note that the above real cylindrical cascade is a topological factor of $\tilde F_a$ and hence, $\tilde F_a$ is far from exhibiting any parabolic behavior.  By the other hand, Besikovitch shown that cylindrical cascades are never minimal. 
\end{enumerate} 
\subsection{Infinitely reducible maps}\label{infinitoreducible}
In this section we look for the possibility of conjugate a parabolic fibred map to the fibred identity map  \[Id_{\alpha}(\theta, z)=(\theta+\alpha, z)\] via a fibred holomorphic change of coordinates. Just from the formal point of view, there should exist $H(\theta, z)=\left(\theta, z+h_2(\theta)z^2+h_3(\theta)z^3+\dots\right)$ such that
\begin{equation}\label{conj33}
F\circ H=H\circ Id_{\alpha}.
\end{equation}
By writing out the formal power series in the above equality we get a recursive definition for the coefficients $h_k$:
\begin{displaymath}
\begin{array}{crcl}
(1)& h_1(\theta)&=&1\\
(2)& h_2(\theta+\alpha)-h_2(\theta)&=&a_2(\theta)\\
\vdots&\vdots &=&\vdots \\
(ec_k)&h_k(\theta+\alpha)-h_k(\theta)&=&\sum_{j=2}^{k}a_j\left\{\sum_{r_1+\dots+r_j=k}h_{r_1}\cdots h_{r_j}\right\}(\theta).
\end{array}
\end{displaymath}
However, each cohomological equation $(ec_k)$ has sense only if equations $(2), \dots , (ec_{k-1})$ have continuous solutions $h_2, \dots, h_{k-1}$. A necessary condition for the existence of the coefficients for $H$ is  the vanishing of every mean
\begin{displaymath}
\begin{array}{crcl}
(Hyp_{k})&\int_{\T} \sum_{j=2}^{k}a_j\left\{\sum_{r_1+\dots+r_j=k}h_{r_1}\cdots h_{r_j}\right\}(\theta)d\theta&=&0.
\end{array}
\end{displaymath}
Of course, the hypothesis $(Hyp_{k})$ has a sense only  provided that $(Hyp_{2}), \dots, (Hyp_{k-1})$ hold and  equations $(2), \dots , (ec_{k-1})$ have continuous solutions $h_2, \dots, h_{k-1}$.  This formal computation suggests an algorithm for studying the dynamics of $F$ near the invariant curve:
\begin{enumerate}
\item If $(Hyp_k)$ does not hold, then we can use the continuous change of coordinates $H^{k-1}(\theta, z)=\left(\theta, z+h_2(\theta)z^2+\dots+h_{k-1}(\theta)z^{k-1}\right)$ in order to conjugate the original map $F$ to 
\[
F_{k}(\theta, z)=\left(\theta+\alpha, z+\tilde a_{k}z^k+\dots\right)
\]
with $\int_{\T}{\tilde a_{k}}\neq 0$. Hence, $F_k$, and a posteriori $F$,  presents a $k-1$ petals parabolic dynamics.
\item If $(Hyp_k)$ holds and there exists a continuous solution $h_k$ for $(ec_k)$, we iterate this algorithm for $(Hyp_{k+1})$ and $(ec_{k+1})$.
\item If $(Hyp_k)$ holds but $(ec_k)$ do not admit a continuos solution, then the dynamics of $F$ can be as {\it strange} as a cylindrical cascade (see Example 5 in Section \ref{ejemplos}) and we can not say much more. We call $F$ a {\it non-reducible} map.

\end{enumerate}
If this algorithm stops for some $k$, then we have a sufficient understood of the local dynamics of $F$ around the invariant curve. In the other hand, it may occur that at each step we fall on the point {\it (ii)} of the above algorithm. In this case we say that $F$ is {\it infinitely reducible.} 

A natural question in the infinitely reducible case is whether or no $F$ is conjugated to the fibred identity map $I_{\alpha}$. At least formally, this is true due to (\ref{conj33}). The next example shows that in order to get a topological conjugacy (that is, the uniform convergence of the series $z+h_{2}(\theta)z^2+\dots$), we need to require additional  hypotheses:

\begin{prop}
Let $\alpha$ be  an irrational number not belonging to the Brjuno class.  Then the fibred map
\[
Q(\theta, z)=\left(\theta+\alpha, z+e^{i(\theta-\alpha)}z^2\right)
\]
is infinitely reducible but not topologically conjugated to the fibred identity.
\end{prop}
\noindent
{\it Proof. }  By  Yoccoz's Theorem, $Q$ is not topologically conjugated to the fibred identity. We need to show that $Q$ is infinitely reducible.   Indeed,  every $(Hyp_k)$ holds, otherwise,  Proposition \ref{flor} implies parabolic behavior. Moreover, each equation $(ec_k)$ has a continuous solution since every $\{h_{j}, a_j\}_{j< k}$ are trigonometric polynomials.  $\quad_{\blacksquare}$\\

Under good hypotheses on $F$ and $\alpha$, the situation in rather simple and we obtain the following dichotomy, which represents the main result on this work:
\begin{theo}\label{identidaddiofantina}
Let 
\begin{eqnarray*}
F:B_{\delta}\times \D&\longrightarrow& B_{\delta}\times \C\\
(\theta, z)&\longmapsto&\left(\theta+\alpha, z+a_2(\theta)z^2+\dots\right)
\end{eqnarray*}
be a fibred holomorphic map, analityc on $B_{\delta}\times \D$, where $B_{\delta}=\{\theta \in \C/\Z \ \big| \ Im(\theta)<\delta\}$ for some $\delta>0$ and $\D=\{z\in \C\ \big| \ |z|<1\}$.  Also suppose that $\alpha$ verifies a Diophantine arithmetic condition $CD(c,\tau)$ for some $c>0, \tau \geq 0$. One of the following statements hols:
\begin{enumerate}
\item $F$ is not infinitely reducible and there exists $k\geq 1$ such that $F:B_{\delta}\times \D\rightarrow B_{\delta}\times \C$ presents a $k$ petals parabolic behavior around the  invariant curve.
\item $F$ is infinitely reducible and there exists 
\begin{eqnarray*}
H:B_{\delta/2}\times \D_r&\longrightarrow& B_{\delta/2}\times \C\\
(\theta, z)&\longmapsto&\left(\theta, z+h_2(\theta)z^2+\dots\right)
\end{eqnarray*}
analityc on $B_{\delta/2}\times \D_r$ for some $0<r\leq 1$, such that 
\[
H^{-1}\circ F\circ H= Id_{\alpha}.
\]
\end{enumerate}
\end{theo}
\noindent{\it Proof. } The first part follows easily from Lemma \ref{babykam}  and Proposition \ref{flor}.  For the proof of the second part we will need to estimate the growth of the coefficients $h_{k}$. For that, we will follow closely the original proof by Siegel of the linearization theorem for holomorphic germs with Diophantine rotation number. We start by recalling   some technical lemmas from  the Siegel's paper \cite{SIEG42}. 
\paragraph{Siegel's Lemmas. }
Let $\{\eps_n\}_{n\in \N}$ be  a real sequence verifying 
\begin{equation}\label{condicioneps}
   \eps_n<(2n)^{\nu}
  \end{equation}
  for some $\nu>0$.   Let's define the sequence $\vartheta_1, \vartheta_2, \dots$ recursively by
$
  \vartheta_1=1.
  $
For $k>1$, denote by $\mu_k$ the biggest value among the products $\vartheta_{r_1}\vartheta_{r_2}\cdots\vartheta_{r_j}$ with 
\[
r_1+r_2+\dots+r_j=k>r_1\geq r_2\geq \dots \geq r_j\geq 1
\]
and $2\leq j\leq k$. Define $\vartheta_k=\eps_{k-1}\mu_k$.
\begin{lema}\label{siegel1}
The following holds 
\[
\vartheta_k\leq k^{-2\nu}2^{(5\nu +1)(k-1)}
\]
for every $k\geq 1\quad_{\blacksquare}$
\end{lema}
 Let's define the sequence $\tau_1, \tau_2, \dots$ recursively by
  \begin{eqnarray*}
  \tau_1&=&1\\
  \tau_k&=&\sum \tau_{r_1}\tau_{r_2}\cdots\tau_{r_j}
  \end{eqnarray*}
  where the above sum is taken over every integer solution of $r_1+r_2+\dots+r_j=k$ with $2\leq j\leq k$.
  \begin{lema}\label{lema2}
  The power series 
  \[
  \sum_{k=1}^{\infty} \tau_kz^k
  \]
  converges on the disc $|z|<3-2\sqrt{2}\quad_{\blacksquare}$
  \end{lema}
  Let's define the sequence $\gamma_1, \gamma_2, \dots$ recursively by
  \begin{eqnarray*}
  \gamma_1&=&1\\
  \gamma_k&=&\eps_{k-1}\sum \gamma_{r_1}\gamma_{r_2}\cdots\gamma_{r_j}
  \end{eqnarray*}
  where the above sum is taken over every integer solution of $r_1+r_2+\dots+r_j=k$ with $2\leq j\leq k$.
  \begin{lema}\label{lema3}
  The following inequality holds
  \[
  \gamma_{k}\leq \vartheta_k\tau_k. 
  \]
  Consequently, the power series 
  \[
  \sum_{k=1}^{\infty}\gamma_kz^k
  \]
  converges on the disc $|z|<(3-2\sqrt{2})2^{-5\nu-1}\quad_{\blacksquare}$

  \end{lema}
  \paragraph{Estimates for the analytic cohomological equation. } Consider the cohomological equation
  \begin{equation}\label{ecuacionanalitica}
  h(\theta+\alpha)-h(\theta)=g(\theta)
  \end{equation}
  where $g:B_{\delta}\to \C$ is analityc on the strip $B_{\delta}$ and $\alpha$ verifies the Diophantine condition $CD(c, \tau)$. Also assume that $\int_{\T}g(\theta)d\theta=0$. The classical method of the Fourier series gives the formal solution
  \begin{equation}\label{solformal}
  h(\theta)=\sum_{n\in \Z\setminus \{0\}}\frac{\hat{g}(n)}{e^{in\alpha}-1}e^{in\theta}
  \end{equation}
  where $\hat{g}(n)$ means for the Fourier coefficient. We pick $\hat h(0)=\int_{\T}h(\theta)d\theta=0$. 
  \begin{lema}\label{lema4}
  Denote by $\|g\|_{\delta}=\sup_{\theta\in B_{\delta}}|g(\theta)|$. Suppose $\|g\|_{\delta}\leq E$. With the above hypothesis and notation one has that the solution $h$ is analytic on the strip $B_{\delta}$. Furthermore,  there exists a constant $C=C(c)$ such that for every $d<\delta$ one has 
  \[
  \|h\|_{\delta-d}\leq \frac{EC}{d^{3+\tau}}.
  \]
  \end{lema}
\noindent{\it Proof. }The Fourier coefficient can be estimated by
\[
|\hat g(n)|\leq Ee^{-|n|\delta}
\]
for $n\in \Z\setminus\{0\}$. The Diophantine condition yields
\[
|e^{in\alpha}-1|\leq \frac{C}{|n|^{2+\tau}}
\]
for some constant $C=C(c)$. Using (\ref{solformal}) and the above estimates one gets
\[
\|h\|_{{\delta-d}}\leq EC\sum_{n\in \Z\setminus\{0\}}e^{-|n|d}|n|^{2+\tau}.
\]
Lemma \ref{lema5} bellow completes the proof $\quad_{\blacksquare}$
\begin{lema}\label{lema5}
 For every $s>0$ there exists $C=C(s)$ such that for every  $x\in (0,1)$ one has
 \[
 \sum_{n\geq 0}x^nn^s\leq \frac{C}{(1-x)^{s+1}}
 \]
\end{lema}
 \noindent
 {\it Proof. } Note that
 \begin{eqnarray*}
 \sum_{n\geq 0}x^nn^s&<&\sum_{n\geq 0}(n+s)(n+s-1)\cdots(n+1)x^n\\
 &=&\frac{\partial^s}{\partial x^s}\left(\sum_{n\geq 0}x^{n+s}\right)\\
 &=&\frac{\partial^s}{\partial x^s}\left(\frac{x^s}{1-x}\right)\quad_{\blacksquare}
 \end{eqnarray*}
Now, we can give the proof of the Theorem \ref{identidaddiofantina}. Let $\nu=\nu(\delta, c, \tau)$ be such that
\begin{equation}\label{sumadelta}
\sum_{k\geq 2}\left(\frac{C(k-1)}{(2(k-1))^{\nu}}\right)^{\frac{1}{3+\tau}}<\frac{\delta}{2}.
\end{equation}
For $k\geq 2$ define 
\[
d_k=\left(\frac{C(k-1)}{(2(k-1))^{\nu}}\right)^{\frac{1}{3+\tau}}\quad\textrm{and}\quad \eps_{k-1}=\frac{C(k-1)}{d_k^{3+\tau}}.
\]
By construction, condition (\ref{condicioneps}) holds. Put $\delta_1=\delta$ and recursively $\delta_{k}=\delta_{k-1}-d_k$. Note that by (\ref{sumadelta}) we have $\delta_k>\frac{\delta}{2}$ for every $k\geq 1$.  As $F(\theta, \cdot)$ has an uniform convergence radius, there exists $a>0$ such that $|a_j(\theta)|<a^{j-1}$ for every $j\geq 2$ and $\theta \in B_{\delta}$. By considering the coordinates change $z\mapsto \frac{z}{a}$ we can assume that $|a_j(\theta)|\leq 1$ for every $j\geq 2, \theta \in B_{\delta}$. 
\begin{lema}\label{tamanodeh}
For every $k\geq 1$ the following holds
\[
\|h_k\|_{\delta_k}\leq \gamma_k.
\]
\end{lema}
\noindent
{\it Proof. } The desired inequality holds for $k=1$. Assume the result for every $j<k$. Lemma \ref{lema4} asserts that the solution $h_k$ for the equation $(ec_k)$ verifies
\begin{eqnarray*}
\|h_k\|_{\delta_{k-1}-d_k}&\leq& \frac{C}{d_k^{3+\tau}}\left\|\sum_{j=2}^{k}a_j\left\{\sum_{r_1+\dots+r_j=k}h_{r_1}\cdots h_{r_j}\right\}\right\|_{\delta_{k-1}}\\
\|h_k\|_{\delta_k}&\leq& \frac{C(k-1)}{d_k^{3+\tau}}\sum_{r_1+\dots+r_j=k}\gamma_{r_1}\cdots \gamma_{r_j}\\
&=&\eps_{k-1}\sum_{r_1+\dots+r_j=k}\gamma_{r_1}\cdots \gamma_{r_j}=\gamma_k\quad_{\blacksquare}
\end{eqnarray*}
Finally, putting together Lemmas \ref{lema3} and \ref{tamanodeh}  the proof of the Theorem \ref{identidaddiofantina} is complete $\quad_{\blacksquare}$
\bibliographystyle{plain}
  \bibliography{fibred_petals2.bib}

\end{document}